\documentclass[12pt,a4paper]{article}
\usepackage[margin=15mm]{geometry}
\usepackage[T1]{fontenc}
\usepackage{lmodern}
\usepackage{amsmath,amssymb,amsthm}
\usepackage{microtype}
\usepackage{tikz}
\usepackage{float}
\usepackage{needspace}
\usepackage[hidelinks]{hyperref}
\hypersetup{pdftitle={A proof of the equitable vertex arboricity conjecture}}
\newtheorem{theorem}{Theorem}[section]
\newtheorem{lemma}{Lemma}[section]
\newtheorem{conjecture}{Conjecture}[section]
\newcommand{\A}{\mathcal A}
\newcommand{\T}{\mathcal T}
\newcommand{\OO}{\mathcal O}
\newcommand{\LL}{\mathcal L}
\begin{document}

\title{A proof of the equitable vertex arboricity conjecture\thanks{Supported by the Youth Project for Natural Science Basic Research Program of Shaanxi (No.~2025JC-YBQN-107), and the Fundamental Research Funds for the Central Universities (Nos.~1608723037, 1960524154).}}
\author{Bei Niu\\
Department of Mathematics\\
Xi'an University of Architecture and Technology\\
Xi'an 710055, China\\
\texttt{beiniu93\char64{}163.com}
\and
Xin Zhang\thanks{Corresponding author.}\\
School of Mathematics and Statistics\\
Xidian University\\
Xi'an 710071, China\\
\texttt{xzhang\char64{}xidian.edu.cn}}
\date{}

\maketitle

\begin{abstract}
We prove that every simple graph $G$ admits an equitable partition into $r$ induced forests whenever
$r\ge\lceil(\Delta(G)+1)/2\rceil$. This settles the equitable vertex arboricity conjecture of Wu, Zhang, and Li (2013) in its full form.
\end{abstract}

\section{Introduction}
All graphs are finite and simple. A \emph{tree-$r$-coloring} of a graph $G$ is a partition of $V(G)$ into $r$ classes, each of which induces a forest. Empty classes are allowed. The coloring is \emph{equitable} if the sizes of any two classes differ by at most one. We write $\Delta(G)$ for the maximum degree of $G$.

The \emph{vertex arboricity} of a graph is the least number of classes
in a tree coloring. Chartrand, Kronk, and Wall~\cite{CKW} proved that
every graph $G$ has vertex arboricity at most
$\lceil(\Delta(G)+1)/2\rceil$. Wu, Zhang, and Li~\cite{WZL}
conjectured that this bound also suffices when a tree coloring is
required to be equitable.

\begin{conjecture}\label{conj:EVAC}
Every graph $G$ has an equitable tree-$r$-coloring for every integer
$r\ge\lceil(\Delta(G)+1)/2\rceil$.
\end{conjecture}

In this paper, we prove the conjecture in the following equivalent form.

\begin{theorem}\label{thm:main}
Let $r$ be a positive integer. Every graph $G$ with $\Delta(G)\le 2r-1$
admits an equitable tree-$r$-coloring.
\end{theorem}

\noindent The degree bound is sharp, since $K_{2r+1}$ has maximum degree $2r$
and every induced forest in it has at most two vertices.

\subsection*{Related work}

Several earlier results establish Conjecture~\ref{conj:EVAC} under
restrictions on the maximum degree. Zhang~\cite{Zhang} proved it for
graphs of maximum degree at most three, and Chen, Liu, and
Zhang~\cite{CLZ} extended the result to maximum degree at most six.
Niu and Zhang~\cite{NZ} obtained equitable tree-$4$-colorings for
graphs of maximum degree at most seven and equitable tree-$5$-colorings
for graphs of maximum degree at most nine.

Other results allow arbitrary maximum degree but impose a condition
on the order. Zhang and Wu~\cite{ZW} proved that every graph $G$ with
$\Delta(G)\ge |V(G)|/2$ has an equitable tree-coloring with
$\lceil(\Delta(G)+1)/2\rceil$ colors. Zhang and Niu~\cite{ZN}
proved that an $n$-vertex graph of maximum degree at most $2r-1$
has an equitable partition into $r$ induced linear forests when
$n\le4r$. 
Gao, Lo, and Shan~\cite{GLS} proved that if $\Delta(G)\le D$ and
$|V(G)|\ge3D^4$, then every integer $r\ge(D+2)/2$ suffices for an
equitable tree-$r$-coloring. For even $D$, their bound gives the
conjectured threshold under the order hypothesis. For $D=2k-1$, it
gives $k+1$ colors, while the conjecture asks for $k$.
Our Theorem~\ref{thm:main} gives the conjectured bound for every graph,
with no additional restriction on its maximum degree or order.

The conjecture has also been verified for several classes of graphs.
Chen, Gao, Shan, Wang, and Wu~\cite{CGSWW} proved it for all
$5$-degenerate graphs. Zhang, Niu, Li, and
Li~\cite{ZNLL} proved it for $d$-degenerate graphs with
$\Delta(G)\ge9.818d$, and Niu, Li, and Zhang~\cite{NLZ} proved it
for interval graphs. Esperet, Lemoine, and Maffray~\cite{ELM} also
obtained a bound depending only on the degeneracy.
These results are related to the equitable partitions into subgraphs
of lower degeneracy studied by Kostochka, Nakprasit, and
Pemmaraju~\cite{KNP}.

For planar graphs, Esperet, Lemoine, and Maffray~\cite{ELM}
proved that every $r\ge4$ suffices for an equitable
tree-$r$-coloring, independently of the maximum degree.
Further results for outerplanar graphs and planar graphs with
restrictions on short cycles appear in~\cite{TL,WZL,Zplanar}.
Kim, Oum, and Zhang~\cite{KOZ} studied related equitable
partitions of planar graphs into subgraphs of low degeneracy
and showed that every planar graph has an equitable partition
into three parts, two of which induce forests.
Equitable tree-colorings have also been studied for plane graphs
with independent crossings~\cite{NZG} and squares of
cycles~\cite{LiuZhang}.

Further work concerns equitable tree-colorings in which the induced
forests have bounded maximum degree, and the least threshold above
which every number of colors is possible. These parameters have been
studied for complete bipartite graphs~\cite{MGZY,NN2,TJT,TL,WZL},
complete tripartite graphs~\cite{GZM2,GZM3,NN2}, balanced complete
multipartite graphs~\cite{LNNR}, and Cartesian and corona
products~\cite{GMJL,Niu2025}. General bounds and extremal questions
were considered by Mao, Guo, Lai, and Zhao~\cite{MGLZ} and
Guo~\cite{GuoNG}. When the induced forests have maximum degree at
most one, the problem is an instance of equitable defective coloring.
Related results appear in~\cite{FKLMWZ,GZ,LiZhang,WVY}.

Besides the existence of equitable tree-colorings, researchers have
studied the complexity of deciding whether a graph admits such a
coloring with a prescribed number of colors.
Nakprasit and Nakprasit~\cite{NNcomplexity} proved that this
problem is NP-complete for general graphs and gave polynomial-time
criteria for complete bipartite graphs. Li and Zhang~\cite{LiZhangTW}
studied the same problem on graphs of bounded treewidth, obtaining
both algorithms and parameterized hardness results.
Niu, Li, and Zhang~\cite{NLZ} proved hardness results for restricted
classes of interval graphs and gave a linear-time algorithm for
deciding whether a proper interval graph admits an equitable
tree-coloring with a prescribed number of colors. They also gave
a $2$-approximation algorithm for minimizing the number of colors
in an equitable tree-coloring of a chordal graph.

Algorithmic results are also available for variants of this problem.
For equitable tree-colorings in which each color class induces a
forest of bounded maximum degree, algorithms for complete bipartite
and tripartite graphs appear in~\cite{NN2}. Niu~\cite{Niu2025}
studied equitable tree-coloring parameters and algorithms for corona
products. More generally, Zhang, Zhang, Niu, and Li~\cite{ZZNL}
gave algorithms for partitioning degenerate graphs into subgraphs
of lower degeneracy.

There is also a list-coloring version, introduced by
Zhang~\cite{Zlist}, in which each vertex has a prescribed list of
colors and each color is used at most $\lceil |V(G)|/k\rceil$ times
for lists of size $k$. Equitable list arboricity has been studied for
grids~\cite{DDFS}, graphs of bounded treewidth~\cite{LiZhangList},
and powers of paths and cycles~\cite{KMP}. More general equitable
list partitions into subgraphs of bounded degeneracy were studied
in~\cite{DFS,ZZmodules}.

For comparison, the well-known Hajnal--Szemer\'edi theorem~\cite{HS} guarantees
an equitable proper coloring with any number of colors at least
$\Delta(G)+1$. Shorter and algorithmic proofs appear
in~\cite{KK,KKMS}. Kierstead and Kostochka~\cite{KKOre} extended
this to the Ore-type condition $d_G(u)+d_G(v)\le2r-1$ for every
edge $uv$, which guarantees an equitable proper $r$-coloring.
The possibility of using $\Delta(G)$ colors was studied by Chen,
Lih, and Wu~\cite{CLW}; further background on equitable proper
coloring can be found in~\cite{KKX,KKY}.

\vspace{2mm}\noindent \textbf{Organization:} Section~\ref{sec:extremal} develops the recoloring tools used in the
proof and bounds the number of color classes involved.
Section~\ref{sec:weighted} uses a degree count to find a configuration
needed for a further recoloring. Section~\ref{sec:recolor} carries out
this recoloring and derives restrictions on the color classes.
Section~\ref{sec:final} uses these restrictions in a final degree count
to complete the proof.

\vspace{2mm}\noindent  \textbf{Notations}: For a vertex set $S$, write $G[S]$ for the induced subgraph and $d_S(v)=|N_G(v)\cap S|$. We use $d_G(v)$ for the degree in $G$ and $d_H(v)$ for the degree in a specified subgraph $H$. For disjoint sets $S,T$, let $e(S,T)$ denote the number of edges between them.

\section{A coloring with one uncolored vertex}\label{sec:extremal}
We prove Theorem~\ref{thm:main} by induction on $r$. If $r=1$, a graph of maximum degree at most one is a forest. Fix $r\ge2$ and assume the theorem for every smaller positive integer. Suppose it fails for $r$, and choose a counterexample $G$ of minimum order. Put $c=\lfloor(|V(G)|-1)/r\rfloor$. For any vertex $z$, the graph $G-z$ has an equitable tree-$r$-coloring, with class sizes $c$ or $c+1$. We call this a \emph{partial coloring} of $G$ with $z$ uncolored. It suffices to find a class of size $c$ into which $z$ can be inserted without creating a cycle.

A vertex $v$ outside a class $V$ is \emph{movable to $V$} if $G[V\cup\{v\}]$ is a forest. This holds exactly when each component of $G[V]$ contains at most one neighbor of $v$. Indeed, two neighbors in one component are joined by a path that forms a cycle with $v$, and every new cycle arises this way. Write $V^*$ for the set of nonisolated vertices of $G[V]$, so that
\begin{equation}\label{eq:nonisolated}
V^*=\{w\in V:d_V(w)\ge1\}.
\end{equation} If $v$ is not movable to $V$, two of its neighbors lie in one nontrivial component and therefore belong to $V^*$.

The \emph{auxiliary digraph} $D$ has the color classes as vertices. For distinct classes $V,Z$, it has an arc $V\to Z$ if a vertex of $V$ is movable to $Z$. Such a vertex is a \emph{witness} of the arc.
To \emph{switch along a directed path}, choose a witness for each arc
$V\to Z$ and simultaneously move it from $V$ to $Z$. If the path has
positive length, the first class loses one vertex, the last gains one,
and every internal class loses one and gains one. All resulting classes
induce forests. Indeed, the first class only loses a vertex, and every
other class receives a vertex movable to the original class and may
also lose a vertex. A switch along a path of length zero makes no
changes.

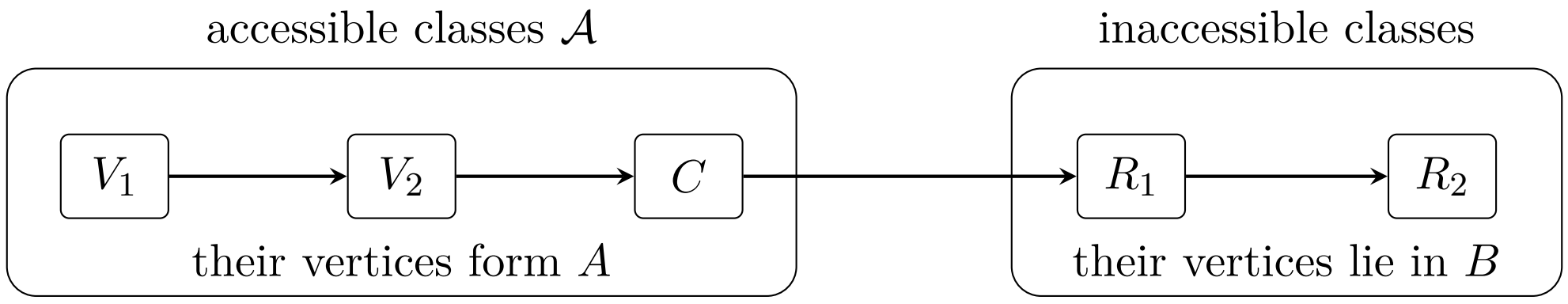
\begin{figure}[H]
\centering
\begin{tikzpicture}[x=1cm,y=1cm,font=\small,>=stealth,
  class/.style={draw,rounded corners=2pt,fill=white,
    minimum width=9mm,minimum height=7mm}]
  \draw[rounded corners=7pt] (-.7,-1) rectangle (5.9,.9);
  \draw[rounded corners=7pt] (7.7,-1) rectangle (12.3,.9);
  \node at (2.6,1.25) {accessible classes $\A$};
  \node at (10,1.25) {inaccessible classes};
  \node[class] (V1) at (.2,0) {$V_1$};
  \node[class] (V2) at (2.6,0) {$V_2$};
  \node[class] (C) at (5,0) {$C$};
  \node[class] (R1) at (8.7,0) {$R_1$};
  \node[class] (R2) at (11.3,0) {$R_2$};
  \draw[->,thick] (V1)--(V2);
  \draw[->,thick] (V2)--(C);
  \draw[->,thick] (C)--(R1);
  \draw[->,thick] (R1)--(R2);
  \node at (2.6,-.7) {their vertices form $A$};
  \node at (10,-.7) {their vertices lie in $B$};
\end{tikzpicture}
\caption{An example of accessibility in $D$.
Each box represents a color class, and arrows are arcs of $D$.
The class $V_1$ reaches $C$ through $V_2$, while $C$ reaches itself
by a path of length zero. Neither $R_1$ nor $R_2$ has a directed
path to $C$. The arc $C\to R_1$ does not make $R_1$ accessible.
No arc can go from an inaccessible class to an accessible class,
since it would give a path to $C$.
The uncolored vertex belongs to $B$ but is not a vertex of $D$
and is not drawn.}
\label{fig:accessibleclasses}
\end{figure}

In a partial coloring, choose a class $C$ of size $c$. A class is \emph{accessible} if it has a directed path to $C$, allowing the path of length zero. Write $\A$ for the family of accessible classes and put
\begin{equation}\label{eq:partition}
a=|\A|,\qquad b=r-a,\qquad A=\bigcup_{V\in\A}V,\qquad B=V(G)\setminus A.
\end{equation}
Thus $B$ contains the uncolored vertex and all vertices of the $b$ inaccessible classes.
Figure~\ref{fig:accessibleclasses} illustrates accessibility in the
auxiliary digraph.

Among all partial colorings and choices of $C$, maximize $a$ and then minimize $\sum_{V\in\A}|V^*|$. We keep the class labels when their vertex sets change.

\begin{lemma}\label{lem:absorb}
We have $a\ge2$ and $1\le b<r$. For every $V\in\A$ and $x\in B$,
\begin{equation}\label{eq:basic}
c\le |V|\le c+1,\qquad |B|=b(c+1)+1\ge b|V|+1,\qquad
d_{V^*}(x)\ge2,\qquad d_B(x)\le2b-1.
\end{equation}
For each $x\in B$, the set $B\setminus\{x\}$ can be partitioned into
$b$ induced forests of size $c+1$. Consequently, any partition of
$A\cup\{x\}$ into $a$ induced forests of sizes $c$ or $c+1$ extends to
an equitable tree-$r$-coloring of $G$.
\end{lemma}
\begin{proof}
No colored vertex in $B$ is movable to an accessible class $V$, since its own class would then reach $C$. The uncolored vertex cannot enter $V$ either. If it could, we could switch along a path from $V$ to $C$ and insert it in $V$. Only $C$ would gain a vertex, giving an equitable coloring. This also covers $V=C$, when no switch is needed.

Thus each $x\in B$ has at least two neighbors in $V^*$ for every $V\in\A$, and $d_B(x)\le2r-1-2a=2b-1$. The set $B$ contains the uncolored vertex, so $d_B(x)\ge0$ for some
$x\in B$. The bound $d_B(x)\le2b-1$ therefore gives $b\ge1$. Also,
$C\in\A$ gives $a\ge1$, and hence $b=r-a<r$.

If any inaccessible class had size $c$, the $b$ inaccessible classes
together with the uncolored vertex would give
$bc<|B|\le b(c+1)$. Since $1\le b<r$ and
$\Delta(G[B])\le2b-1$, the induction hypothesis gives an equitable
tree-$b$-coloring of $G[B]$. The bounds on $|B|$ imply that each of
its classes has size $c$ or $c+1$. Together with the accessible
classes, they give an equitable tree-$r$-coloring of $G$.
Hence each inaccessible class has size $c+1$ and $|B|=b(c+1)+1$.

For every $x\in B$, the graph $G[B\setminus\{x\}]$ has
$b(c+1)$ vertices and maximum degree at most $2b-1$. By the
induction hypothesis, it has an equitable tree-$b$-coloring, in which
every class has size $c+1$. Thus any partition of $A\cup\{x\}$ into
$a$ induced forests of sizes $c$ or $c+1$ completes an equitable
tree-$r$-coloring of $G$.

Finally, we show $a\ge 2$. If $a=1$, then $C$ is the sole accessible class and $b=r-1$. By~\eqref{eq:nonisolated}, each vertex of $C^*$ has a neighbor in $C$, so it has at most $\Delta(G)-1\le 2r-2=2b$ neighbors in $B$. Counting edges at both ends gives
\[
2|B|\le e(B,C^*)\le 2b|C^*|,
\]
contrary to $|B|=b(c+1)+1>b|C^*|$.
\end{proof}

For each $U\in\A$, let $\T_U$ be the family of classes in
$\A\setminus\{U\}$ that cannot reach $C$ in $D[\A]-U$. When $U=C$,
set $\T_C=\A\setminus\{C\}$. Since $a\ge2$ by Lemma~\ref{lem:absorb}, the family $\T_C$ is
nonempty. Among the finitely many nonempty families $\T_U$ with $U\in\A$,
choose one of minimum size.
Put
\begin{equation}\label{eq:families}
\T=\T_U,\qquad t=|\T|,\qquad
\OO=\A\setminus(\T\cup\{U\}).
\end{equation}
Then $1\le t\le a-1$, and neither $U$ nor $C$ belongs to $\T$.

\begin{lemma}\label{lem:paths}
No arc goes from $\T$ to $\OO$. For every $W\in\T$, each class in
$\A\setminus\{W\}$ reaches $C$ in $D[\A]-W$.
\end{lemma}

\begin{proof}
Every class in $\OO$ has a path to $C$ that avoids $U$. Such a path
also avoids $\T$, because otherwise a class in $\T$ would have a path
to $C$ avoiding $U$. Hence there is no arc from $\T$ to $\OO$.

A path from $U$ to $C$ also avoids $\T$. If it contained a class in
$\T$, the part of the path from that class to $C$ would avoid $U$,
contrary to the definition of $\T$ in~\eqref{eq:families}. When $U=C$, this path has length
zero and avoids $\T$ as well.

Fix $W\in\T$. The paths found above show that every class in
$\OO\cup\{U\}$ reaches $C$ in $D[\A]-W$, since these paths avoid
$\T$ and hence avoid $W$. Therefore
$\T_W\subseteq\T\setminus\{W\}$, which gives $|\T_W|\le t-1$.
If $\T_W$ were nonempty, the choice of $U$ would give
$|\T_W|\ge|\T_U|=t$. Thus $\T_W$ is empty, and every class in
$\A\setminus\{W\}$ reaches $C$ in $D[\A]-W$.
\end{proof}

\begin{figure}[H]
\centering
\begin{tikzpicture}[x=1cm,y=1cm,font=\small,>=stealth,
  class/.style={draw,rounded corners=2pt,fill=white,
    minimum width=8mm,minimum height=7mm}]
  \node[anchor=west] at (-.5,1.9) {(a) $U\ne C$.};
  \draw (1.5,0) ellipse (2cm and 1.2cm);
  \draw (10.5,0) ellipse (2cm and 1.2cm);
  \node at (1.5,.9) {$\T$};
  \node at (10.5,.9) {$\OO$};
  \node[class] (V) at (1.5,0) {$V$};
  \node[class] (U) at (6,0) {$U$};
  \node[class] (C) at (10.5,0) {$C$};
  \draw[->,thick,dashed] (V)--(U);
  \draw[->,thick,dashed] (U)--(C);
  \draw[->] (3,-.8)--(9,-.8);
  \draw[thick] (5.88,-.92)--(6.12,-.68)
    (5.88,-.68)--(6.12,-.92);
  \node[below] at (6,-1) {no arc from $\T$ to $\OO$};
  \begin{scope}[yshift=-3.4cm]
    \node[anchor=west] at (-.5,1.05) {(b) $U=C$.};
    \draw (1.5,0) ellipse (2cm and .75cm);
    \node at (1.5,.5) {$\T$};
    \node[class] (Vb) at (1.5,-.1) {$V$};
    \node[class] (Cb) at (6,-.1) {$U=C$};
    \draw[->,thick,dashed] (Vb)--(Cb);
    \node at (10.5,-.1) {$\OO=\varnothing$};
  \end{scope}
\end{tikzpicture}
\caption{The families in Lemma~\ref{lem:paths}.
Boxes represent classes, and dashed arrows represent directed
paths with intermediate classes omitted. Every path from a class
$V\in\T$ to $C$ passes through $U$. In (a), its classes before
$U$ lie in $\T$, and those after $U$ lie in $\OO$.
The crossed arrow indicates that no arc goes from $\T$ to $\OO$.
In (b), $\T=\A\setminus\{C\}$ and $\OO$ is empty.}
\label{fig:accessibility}
\end{figure}
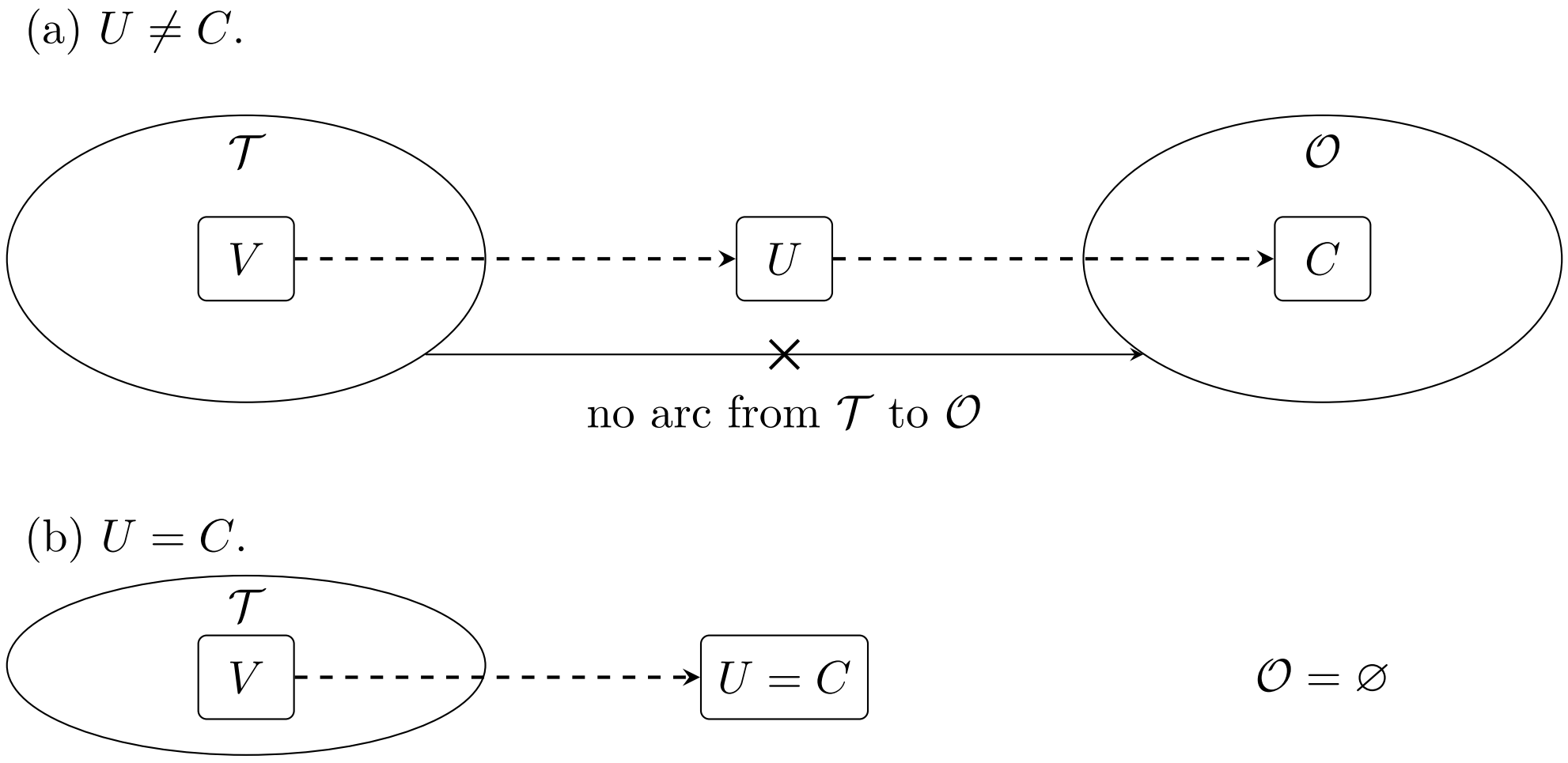

Figure~\ref{fig:accessibility} shows the relation between $\T$, $U$,
and $\OO$ established above.

We next consider replacing a vertex of a class in $\T$ by a vertex
of $B$.

\begin{lemma}\label{lem:exchange}
Let $W\in\T$, $x\in B$, and $w\in W$. If $\widehat W=(W\setminus\{w\})\cup\{x\}$ induces a forest, then $w$ is movable to no class in $\A\setminus\{W\}$, and no vertex of $B\setminus\{x\}$ is movable to $\widehat W$.
\end{lemma}
\begin{proof}
Use Lemma~\ref{lem:absorb} to color $B\setminus\{x\}$ in $b$ classes of size $c+1$. Replace $w$ by $x$ in $W$. We now have a partial coloring with $w$ uncolored and the same class sizes.
By Lemma~\ref{lem:paths}, every class in $\A\setminus\{W\}$ has a
path to $C$ that avoids $W$. These classes are unchanged by the
recoloring, so the same paths remain in the new auxiliary digraph.
If $w$ were movable to any of these classes, we could switch along its path to $C$ and insert $w$, completing the coloring.

Since $W$ is accessible, it has a path to $C$ in the original
auxiliary digraph. The witness for the first arc of this path is
different from $w$, because $w$ is movable to no class in
$\A\setminus\{W\}$. This witness remains in $\widehat W$, and the
other classes on the path are unchanged. Hence $\widehat W$ also
reaches $C$ in the new auxiliary digraph. If a vertex of
$B\setminus\{x\}$ were movable to $\widehat W$, its class would also
reach $C$. This would give more than $a$ accessible classes,
contrary to the choice of the original coloring.
\end{proof}

\needspace{6\baselineskip}
\begin{lemma}\label{lem:tbound}
We have $t<2b$.
\end{lemma}
\begin{proof}
Suppose we can select $2b$ classes from $\T$. For each $x\in B$,
let $j_x$ be the number of selected classes containing exactly two
neighbors of $x$. By Lemma~\ref{lem:absorb}, every accessible class
contains at least two neighbors of $x$. Each of the other $2b-j_x$
selected classes therefore contains at least three. Hence
\[
2(a-2b)+2j_x+3(2b-j_x)+d_B(x)
\le d_G(x)\le2r-1.
\]
Since $a+b=r$ by~\eqref{eq:partition}, this gives
\begin{equation}\label{eq:marks}
j_x\ge1+d_B(x).
\end{equation}

Choose a maximal independent set $X$ in $G[B]$. Every vertex of
$B\setminus X$ has a neighbor in $X$, since otherwise it could be
added to $X$. Thus the closed neighborhoods of the vertices in $X$
cover $B$, and
\begin{equation}\label{eq:cover}
\sum_{x\in X}(1+d_B(x))\ge|B|.
\end{equation}

For each $x\in X$, mark its two neighbors in each of the $j_x$
selected classes. Counting marks with multiplicity, their total number
is $2\sum_{x\in X}j_x$. Since $j_x\ge1+d_B(x)$ by~\eqref{eq:marks} and
$|B|=b(c+1)+1$ by \eqref{eq:basic}, \eqref{eq:cover} gives
\[
2\sum_{x\in X}j_x
\ge2\sum_{x\in X}(1+d_B(x))
\ge2|B|>2b(c+1).
\]
The selected classes contain at most $2b(c+1)$ vertices altogether.
Thus two distinct vertices $x,y\in X$ mark the same vertex $v$ in
some selected class $W$. Each of $x$ and $y$ has exactly two neighbors
in $W$, including $v$.

Replacing $v$ by $x$ leaves a forest, since $x$ has just one
remaining neighbor in $W$. The vertices $x$ and $y$ are nonadjacent,
so $y$ also has just one neighbor in the new class and is movable to
it. This contradicts Lemma~\ref{lem:exchange}.
\end{proof}

\section{Counting edges at an accessible class}\label{sec:weighted}

To find a class $W\in\T$ such that every vertex of $W^*$ has a
neighbor in every other class of $\T$, we define a digraph $D_0$
on $\A$. 
For distinct classes $V,Z\in\A$,
put an arc $V\to Z$ in $D_0$ if some vertex of $V^*$ has no neighbor
in $Z$. This defines a digraph $D_0$ on $\A$.

\begin{lemma}\label{lem:acyclic}
The digraph $D_0$ is acyclic.
\end{lemma}

\begin{proof}
Suppose $D_0$ contains a directed cycle
$V_1\to V_2\to\cdots\to V_m\to V_1$, with indices taken modulo $m$.
For each $i$, choose a witness $x_i\in V_i^*$ with no neighbor in
$V_{i+1}$. Simultaneously move each $x_i$ from $V_i$ to $V_{i+1}$.
For every original class $V$, write $V'$ for the vertex set of the
class with the same label after the switch. Thus
\begin{equation}\label{eq:cycleclasses}
V_i'=(V_i\setminus\{x_i\})\cup\{x_{i-1}\},
\end{equation}
and $R'=R$ for each class $R$ outside the cycle. Each $V_i$ loses
the nonisolated vertex $x_i$ and receives $x_{i-1}$, which has no
neighbor in $V_i$. Hence $V_i'$ induces a forest, $|V_i'|=|V_i|$, and $V_i'$ has
strictly fewer nonisolated vertices than $V_i$.

\medskip
\noindent\textbf{Claim.}
For every $S\in\A$, the class $S'$ reaches $C'$ in
the new auxiliary digraph.

\smallskip
\noindent
First consider the classes on the cycle. The vertex $x_i$ now
belongs to $V_{i+1}'$ and is movable to $V_i'$. Indeed,
\[
V_i'\cup\{x_i\}=V_i\cup\{x_{i-1}\},
\]
which induces a forest since $x_{i-1}$ has no neighbor in $V_i$.
Thus the new auxiliary digraph contains the arc $V_{i+1}'\to V_i'$
for every $i$. These arcs form a reversed cycle, so each $V_i'$ can
reach every other class on this cycle.

Next let $R$ be a class outside the cycle. Suppose the original
digraph contains an arc $R\to V_i$ with witness $z\in R$. Then
$V_i\cup\{z\}$ induces a forest. To obtain $V_i'\cup\{z\}$, we
delete $x_i$ and add $x_{i-1}$ by~\eqref{eq:cycleclasses}. The only possible neighbor of
$x_{i-1}$ in this set is $z$, so $V_i'\cup\{z\}$ also induces a
forest. Hence the new digraph contains the arc $R\to V_i'$.

Now suppose the original digraph contains an arc $V_i\to R$ with
witness $w\in V_i$. Since $R'=R$, the vertex $w$ is still movable
to $R$. If $w=x_i$, then $w\in V_{i+1}'$ and gives the arc
$V_{i+1}'\to R$. If $w\ne x_i$, then $w\in V_i'$ and gives the
arc $V_i'\to R$. Arcs between two classes outside the cycle also
remain, since both classes keep their original vertex sets.

Fix $S\in\A$ and take an original path from $S$ to $C$. If the path
avoids the cycle, it remains a path from $S'=S$ to $C'=C$.
Otherwise, let $V_i$ be the first class on the cycle that the path
meets. Then $S'$ can reach $V_i'$. This is immediate if $S=V_i$.
Otherwise, let $R$ be the class immediately before $V_i$ on the path.
The path from $S$ to $R$ lies outside the cycle, so it remains in the
new digraph. The arc $R\to V_i'$ then gives a path from $S$ to $V_i'$. If $C$ lies on
the cycle, the reversed cycle leads from $V_i'$ to $C'$.
If $C$ lies outside the cycle, let $V_j$ be the last class on the
original path that belongs to the cycle, and let $Z$ be the next
class on the path. The witness for $V_j\to Z$ now belongs to either
$V_j'$ or $V_{j+1}'$ and is still movable to $Z$. Since $V_i'$ can
reach both classes along the reversed cycle, it can reach $Z$.
The rest of the path from $Z$ to $C$ lies outside the cycle, so it
remains in the new digraph. Thus $S'$ reaches $C'=C$, proving the claim.

The class $C'$ still has size $c$. By the claim, all $a$ classes
$S'$ with $S\in\A$ are accessible in the new coloring. By the
maximality of $a$, these are exactly its accessible classes. Each
$V_i$ belongs to $\A$, and each $V_i'$ has fewer nonisolated
vertices than $V_i$. Every accessible class outside the cycle has
the same vertex set as before. Thus the switch strictly decreases
the total number of nonisolated vertices in the accessible classes,
contrary to the choice of the partial coloring.
\end{proof}

By Lemma~\ref{lem:acyclic}, $D_0[\T]$ is acyclic. Since $\T$ is
finite and nonempty, $D_0[\T]$ has a class $W$ of outdegree zero. Thus
\begin{equation}\label{eq:Wneighbors}
d_V(w)\ge1 \qquad (w\in W^*,\ V\in\T\setminus\{W\}).
\end{equation} Indeed, if some $x\in W^*$ had no neighbor
in such a class $V$, then $x$ would witness the arc $W\to V$ in
$D_0[\T]$, contradicting the choice of $W$.

Let $P$ be the set of vertices in $W^*$ that are movable to no other
accessible class, and put
\begin{equation}\label{eq:PQ}
F=G[W^*],\qquad Q=W^*\setminus P.
\end{equation}

Fix $x\in B$. Since $x$ is not movable to $W$ and $G[W]$ is a
forest, $G[W\cup\{x\}]$ contains a cycle through $x$. Removing $x$
from this cycle leaves a path in $G[W]$ joining two distinct
neighbors of $x$. For every $D'\subseteq W$, we have
\begin{equation}\label{eq:pathsmeet}
\begin{aligned}
&G[(W\setminus D')\cup\{x\}]\text{ is a forest}\\
&\qquad\Longleftrightarrow\quad
D'\text{ meets every path in }G[W]\text{ joining two distinct neighbors of }x.
\end{aligned}
\end{equation}
Indeed, every cycle in $G[(W\setminus D')\cup\{x\}]$ contains $x$
and gives a path of this kind avoiding $D'$. Conversely, any such
path avoiding $D'$ forms a cycle with $x$.

Suppose that $w\in W$ lies on every path in $G[W]$ joining two
distinct neighbors of $x$. Since $x$ is not movable to $W$, we may
choose a path $L$ in $G[W]$ joining two distinct neighbors of $x$.
Then $w$ lies on $L$. Since $L$ has positive length, $w$ is incident
with an edge of $L$, so $w\in W^*$.
By~\eqref{eq:pathsmeet} with $D'=\{w\}$, $x$ is movable to
$W\setminus\{w\}$. Lemma~\ref{lem:exchange} therefore implies that
$w$ is movable to no other accessible class. Hence $w\in P$.

A vertex $x\in B$ is \emph{exceptional for $W$} if it has exactly
three neighbors in $W^*$ and all three belong to $Q$.
Following Niu and Zhang~\cite{NZ}, we consider the paths between
three neighbors of $x$. Suppose that $x$ has exactly three neighbors
$u_1,u_2,u_3$ in $W^*$ and that they lie in the same tree of $F$.
The $u_1$--$u_2$, $u_2$--$u_3$, and $u_3$--$u_1$ paths in this
tree have a unique common vertex $m$, called their \emph{median}.
The three paths from $m$ to $u_1$, $u_2$, and $u_3$ are the
\emph{branches} of $x$. Note that the median may be one of $u_1,u_2,u_3$.

\begin{lemma}\label{lem:median}
Let $x\in B$.
\begin{enumerate}
\item If $x$ has exactly two neighbors in $W^*$, then $F$ contains
a path joining them. Every vertex on this path belongs to $P$.

\item If $x$ is exceptional for $W$, then its three neighbors in
$W^*$ lie in the same tree of $F$. Their median $m$ belongs to $P$
and has degree at least three in $F$. If $w_1$ and $w_2$ lie on
distinct branches of $x$ and neither is $m$, then $x$ is movable
to $W\setminus\{w_1,w_2\}$.
\end{enumerate}
\end{lemma}

\begin{proof}
Since $x$ is not movable to $W$, some path in $G[W]$ joins two
distinct neighbors of $x$. Every such path lies in $F$, since the
vertices of $W\setminus W^*$ are isolated in $G[W]$ by~\eqref{eq:nonisolated}.

Suppose that $x$ has exactly two neighbors $u$ and $v$ in $W^*$.
Since every path in $G[W]$ joining two distinct neighbors of $x$
lies in $F$, its endpoints must be $u$ and $v$. The forest $F$
contains exactly one such path. Let $w$ be any vertex on it.
Deleting $w$ breaks the path, so no two neighbors of $x$ are joined
in $G[W\setminus\{w\}]$. Thus $x$ is movable to $W\setminus\{w\}$.
By Lemma~\ref{lem:exchange}, $w$ is movable to no other accessible
class. Since $w\in W^*$, this gives $w\in P$.

Suppose that $x$ is exceptional for $W$. Since some path in $F$ joins
two neighbors of $x$, we may label its three neighbors in $W^*$ as
$u_1,u_2,u_3$ so that $u_1$ and $u_2$ lie in the same tree of $F$.
If $u_3$ lay in a different tree, the $u_1$--$u_2$ path would be
the only path in $G[W]$ joining two distinct neighbors of $x$.
Deleting $u_1$ would break this path, so $x$ would be movable to
$W\setminus\{u_1\}$. By Lemma~\ref{lem:exchange}, $u_1$ would then
be movable to no other accessible class. Since $u_1\in W^*$, this
would give $u_1\in P$. But $u_1\in Q$ by the definition of an
exceptional vertex. Thus $u_3$ lies in the same tree as $u_1$ and $u_2$.

Their median $m$ lies on every path in $G[W]$ joining two distinct
neighbors of $x$. Hence $x$ is movable to $W\setminus\{m\}$ by
\eqref{eq:pathsmeet}. Lemma~\ref{lem:exchange} then implies that
$m$ is movable to no other accessible class. Since $m\in W^*$,
we obtain $m\in P$.
Since $u_1,u_2,u_3$ belong to $Q$, each is distinct from $m$.
The three branches meet only at $m$ and therefore leave $m$ along
three distinct edges. Thus $m$ has degree at least three in $F$.

Choose vertices $w_1$ and $w_2$ on two distinct branches, with
$w_1,w_2\ne m$. Relabel the neighbors so that $w_1$ lies on the
$m$--$u_1$ path and $w_2$ lies on the $m$--$u_2$ path. Then $w_1$
lies on both the $u_1$--$u_2$ path and the $u_1$--$u_3$ path,
and $w_2$ lies on the $u_2$--$u_3$ path. Thus $\{w_1,w_2\}$
meets every path in $G[W]$ joining two distinct neighbors of $x$.
By~\eqref{eq:pathsmeet}, $x$ is movable to $W\setminus\{w_1,w_2\}$.
\end{proof}

\begin{figure}[H]
\centering
\begin{tikzpicture}[x=1cm,y=1cm,font=\small,
  vertex/.style={circle,fill=black,inner sep=1.7pt},
  removed/.style={circle,draw,fill=white,inner sep=3.2pt}]
  \coordinate (m) at (0,0);
  \coordinate (u1) at (3,1.3);
  \coordinate (u2) at (3,0);
  \coordinate (u3) at (3,-1.3);
  \coordinate (x) at (5,0);
  \coordinate (w1) at (1.5,.65);
  \coordinate (w2) at (1.5,0);
  \draw[thick] (m)--(u1) (m)--(u2) (m)--(u3);
  \draw[dashed] (x)--(u1) (x)--(u2) (x)--(u3);
  \node[vertex,label=left:{$m\in P$}] at (m) {};
  \node[vertex,label=above:{$u_1$}] at (u1) {};
  \node[vertex,label=below right:{$u_2$}] at (u2) {};
  \node[vertex,label=below:{$u_3$}] at (u3) {};
  \node[vertex,label=right:{$x\in B$}] at (x) {};
  \node[removed,label=above:{$w_1$}] at (w1) {};
  \node[removed,label=below:{$w_2$}] at (w2) {};
  \draw[thick] (1.42,.57)--(1.58,.73) (1.42,.73)--(1.58,.57);
  \draw[thick] (1.42,-.08)--(1.58,.08) (1.42,.08)--(1.58,-.08);
\end{tikzpicture}
\caption{The three branches of an exceptional vertex $x$.
Solid lines represent paths in $F$, possibly with additional internal
vertices. Dashed lines are the edges from $x$ to its three neighbors
$u_1,u_2,u_3\in Q$. The crossed vertices $w_1,w_2$ lie on distinct
branches and differ from $m$. Deleting them separates every pair
of these neighbors, allowing $x$ to enter $W$. Only the relevant
paths and edges are shown.}
\label{fig:branches}
\end{figure}
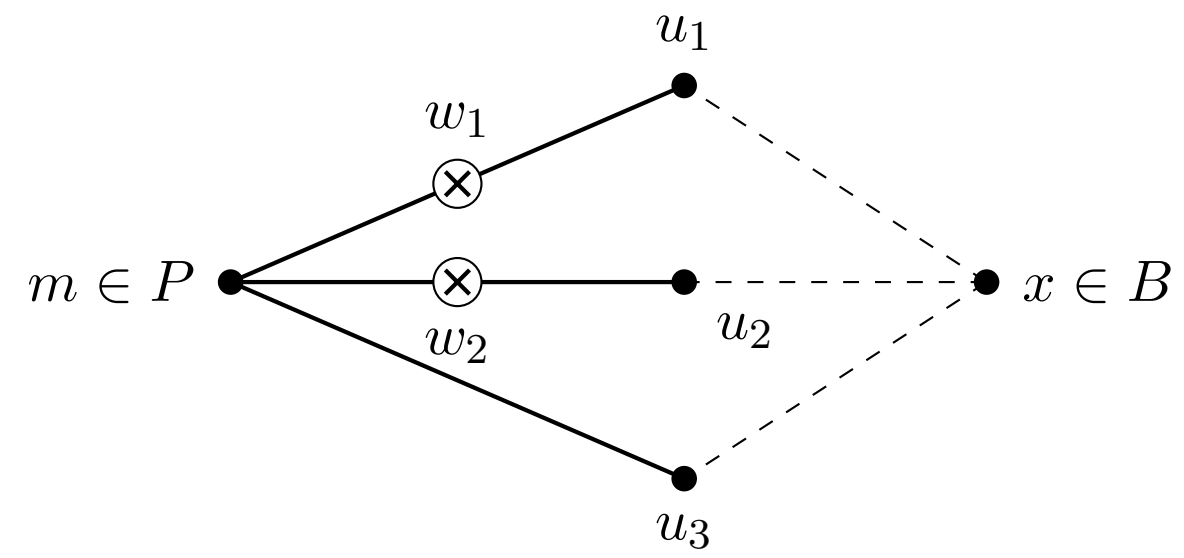

Figure~\ref{fig:branches} illustrates the deletion of two vertices
on distinct branches in Lemma~\ref{lem:median}.

We now count the edges from $B$ to $W^*$. To bound their number at a
vertex $w\in W^*$, we first account for its neighbors in the other
accessible classes. Each $p\in P$ has at least two neighbors in each
such class, since it cannot move there. Recall that $W$ has outdegree zero in $D_0[\T]$ by its choice. If some $q\in Q$
had no neighbor in a class $V\in\T\setminus\{W\}$, then $q\in W^*$
would witness the arc $W\to V$, a contradiction. Thus each $q\in Q$ has a neighbor
in every class of $\T\setminus\{W\}$, as in~\eqref{eq:Wneighbors}. Since $W\in\T$, Lemma~\ref{lem:paths} rules out an arc from $W$ to
any class $V\in\OO$. If $q$ had at most one neighbor in such a class
$V$, then adding $q$ to $V$ could not create a cycle. Thus $q$
would be movable to $V$ and would witness the arc $W\to V$.
Hence $q$ has at least two neighbors in every class of $\OO$.
Since $U\notin\T$ and $U\notin\OO$, for $q\in Q$ we use only the
trivial lower bound $d_U(q)\ge 0$.
The
following table records these lower bounds on $d_V(w)$.
\begin{equation}\label{eq:classtable}
\begin{array}{c|ccc|c}
 & V\in\T\setminus\{W\} & V=U & V\in\OO & \text{sum over }V\\ \hline
w\in P & 2 & 2 & 2 & 2(a-1)\\
w\in Q & 1 & 0 & 2 & (t-1)+2(a-t-1)
\end{array}
\end{equation}
For each $w\in W^*$, let $\sigma(w)$ be the difference between the
actual number of neighbors of $w$ in $A\setminus W$ and the lower
bound in the last column of the table. Thus
\begin{equation}\label{eq:sigmadef}
\sigma(w)=
\begin{cases}
d_{A\setminus W}(w)-2(a-1), & w\in P,\\
d_{A\setminus W}(w)-\bigl((t-1)+2(a-t-1)\bigr), & w\in Q.
\end{cases}
\end{equation}
The bounds in~\eqref{eq:classtable} give
\begin{equation}\label{eq:sigmanonnegative}
\sigma(w)\ge0 \qquad (w\in W^*).
\end{equation}
Put
\begin{equation}\label{eq:sigmatotal}
\sigma_W=\sum_{w\in W^*}\sigma(w),\qquad s=t-b.
\end{equation} By Lemma~\ref{lem:tbound},
\begin{equation}\label{eq:dbound}
s\le b-1.
\end{equation}
For $w\in W^*$, its neighbors in $W$ all belong to $W^*$ by~\eqref{eq:nonisolated}.
Since $F=G[W^*]$ by~\eqref{eq:PQ}, we have
\begin{equation}\label{eq:forestdegrees}
d_F(w)=d_W(w)\ge1,\qquad d_B(w)=d_G(w)-d_F(w)-d_{A\setminus W}(w).
\end{equation}
Using $d_G(w)\le2r-1$,
$a+b=r$ by~\eqref{eq:partition}, and $(t-1)+2(a-t-1)=2a-t-3$, we obtain from~\eqref{eq:sigmadef} and~\eqref{eq:forestdegrees}
\begin{equation}\label{eq:Bdegrees}
\begin{aligned}
d_B(p)&\le 2b+1-d_F(p)-\sigma(p) &&(p\in P),\\
d_B(q)&\le 3b+2+s-d_F(q)-\sigma(q) &&(q\in Q).
\end{aligned}
\end{equation}

For each vertex $x\in B$ exceptional for $W$, let $S_x$ be the set
of all vertices on the three paths from its median to its three
neighbors in $W^*$, including the median and the three neighbors.
Let $J$ be the subgraph of $F$ induced by the union of these sets, so that
\begin{equation}\label{eq:branchsets}
V(J)=\bigcup_{x\text{ exceptional for }W} S_x,\qquad Q_J=Q\cap V(J).
\end{equation}
Let $I$ denote the number of isolated vertices in $G[W]$. Thus
\begin{equation}\label{eq:isolated}
I=|W\setminus W^*|=|W|-|W^*|.
\end{equation}

\begin{lemma}\label{lem:weighted}
We have
\begin{equation}\label{eq:branchbudget}
12+12bI+3\sigma_W
\le \sum_{q\in Q_J}\bigl(4s-\sigma(q)\bigr).
\end{equation}
\end{lemma}

\begin{proof}
Give each edge from $B$ to $W^*$ weight $6$ if its endpoint in
$W^*$ belongs to $P$, weight $4$ if it belongs to $Q_J$, and
weight $3$ if it belongs to $Q\setminus Q_J$.

We first count the weight at each $x\in B$. By~\eqref{eq:basic},
$x$ has at least two neighbors in $W^*$. If it has exactly two,
both belong to $P$ by Lemma~\ref{lem:median}, so the two edges
have total weight $12$. Suppose that $x$ has exactly three neighbors
in $W^*$. If at least one belongs to $P$, the three edges have total
weight at least $6+3+3=12$. Otherwise all three neighbors belong to
$Q$, so $x$ is exceptional for $W$. Each of the three neighbors is
an endpoint of a branch of $x$ and hence belongs to $Q_J$. The three
edges therefore have total weight $4+4+4=12$. Finally, if $x$ has
at least four neighbors in $W^*$, their edges have total weight at
least $4\cdot3=12$. Summing over $B$ gives
\begin{equation}\label{eq:totalweight}
12|B|\le 6e(B,P)+4e(B,Q_J)+3e(B,Q\setminus Q_J).
\end{equation}

We now count the same weight at the vertices of $W^*$. For each
$p\in P$, the $d_B(p)$ edges from $p$ to $B$ have total weight
$6d_B(p)$. Similarly, the edges from $q$ to $B$ have total weight
$4d_B(q)$ if $q\in Q_J$ and $3d_B(q)$ if $q\in Q\setminus Q_J$.
Thus the preceding inequality~\eqref{eq:totalweight} becomes
\begin{equation}\label{eq:vertexweight}
12|B|\le
\sum_{p\in P}6d_B(p)
+\sum_{q\in Q_J}4d_B(q)
+\sum_{q\in Q\setminus Q_J}3d_B(q).
\end{equation}

There is one term on the right of~\eqref{eq:vertexweight} for each $w\in W^*$. We subtract
$12b$ from each term and add $3\sigma(w)$ to it. The total change
is $-12b|W^*|+3\sigma_W$ by~\eqref{eq:sigmatotal}. Making the same change on the left gives
\begin{equation}\label{eq:weightlower}
\begin{aligned}
12|B|-12b|W^*|+3\sigma_W
&\ge 12(b|W|+1)-12b|W^*|+3\sigma_W\\
&=12+12bI+3\sigma_W,
\end{aligned}
\end{equation}
where the inequality follows from~\eqref{eq:basic} and the equality
uses~\eqref{eq:isolated}.

We apply~\eqref{eq:Bdegrees} to the resulting term for each vertex.
For $p\in P$, we obtain
\begin{equation}\label{eq:weightP}
\begin{aligned}
6d_B(p)-12b+3\sigma(p)
&\le 6\bigl(2b+1-d_F(p)-\sigma(p)\bigr)-12b+3\sigma(p)\\
&=-6(d_F(p)-1)-3\sigma(p).
\end{aligned}
\end{equation}
For $q\in Q$, we use $d_B(q)\le 3b+2+s-d_F(q)-\sigma(q)$.
If $q\in Q_J$, its edges to $B$ have weight $4$, giving
\begin{equation}\label{eq:weightQJ}
\begin{aligned}
4d_B(q)-12b+3\sigma(q)
&\le 4\bigl(3b+2+s-d_F(q)-\sigma(q)\bigr)-12b+3\sigma(q)\\
&=4s-\sigma(q)+4(2-d_F(q)).
\end{aligned}
\end{equation}
If $q\in Q\setminus Q_J$, its edges to $B$ have weight $3$, giving
\begin{equation}\label{eq:weightother}
\begin{aligned}
3d_B(q)-12b+3\sigma(q)
&\le 3\bigl(3b+2+s-d_F(q)-\sigma(q)\bigr)-12b+3\sigma(q)\\
&=3(s+2-b-d_F(q))\le0.
\end{aligned}
\end{equation}
The last inequality holds because $s\le b-1$ by~\eqref{eq:dbound} and $d_F(q)\ge1$ by~\eqref{eq:forestdegrees}.

We first sum~\eqref{eq:weightP} over $p\in P$. Since $\sigma(p)\ge0$ by~\eqref{eq:sigmanonnegative},
we may omit $-3\sigma(p)$ for every $p\in P$. For
$p\in P\setminus V(J)$, we may also omit $-6(d_F(p)-1)$, since
$d_F(p)\ge1$ by~\eqref{eq:forestdegrees}. For $p\in P\cap V(J)$,
we retain $-6(d_F(p)-1)$. Thus the sum over $P$ is bounded above by
$-6\sum_{p\in P\cap V(J)}(d_F(p)-1)$.

For $q\in Q_J$, we retain both terms on the right
of~\eqref{eq:weightQJ}. The contribution from each
$q\in Q\setminus Q_J$ is nonpositive by~\eqref{eq:weightother} and
may be omitted. Combining these upper bounds with the lower bound
in~\eqref{eq:weightlower} gives
\begin{equation}\label{eq:degreebalance}
\begin{aligned}
12+12bI+3\sigma_W
&\le \sum_{q\in Q_J}(4s-\sigma(q))\\
&\quad +4\sum_{q\in Q_J}(2-d_F(q))
       -6\sum_{p\in P\cap V(J)}(d_F(p)-1).
\end{aligned}
\end{equation}
To prove~\eqref{eq:branchbudget}, it suffices to show that the last
line of~\eqref{eq:degreebalance} is nonpositive. Since
$Q_J=Q\cap V(J)$ by~\eqref{eq:branchsets}, both sums on that line are over vertices of $J$.
We show that the contribution from each component of $J$ is
nonpositive.

Let $K$ be a component of $J$, and choose $z\in V(K)$. By~\eqref{eq:branchsets}, the vertex $z$ lies on a branch of some
exceptional vertex $x$. The three branches of $x$ lie in $J$ and
meet at their median $m$. Their union is connected and contains
$z$, so all three branches lie in $K$. By Lemma~\ref{lem:median},
$m\in P$. The branches leave $m$ along three distinct edges, so
$d_K(m)\ge3$.

Consider the contribution of $K$ to the last line
of~\eqref{eq:degreebalance}. Since $K$ is a subgraph of $F$, we
have $d_K(w)\le d_F(w)$ for every $w\in V(K)$. Each
$q\in Q\cap V(K)$ belongs to $Q_J$, so its contribution satisfies
\[
4(2-d_F(q))\le 4(2-d_K(q)).
\]
For each $p\in P\cap V(K)$, the corresponding estimate is
\[
-6(d_F(p)-1)
\le -6(d_K(p)-1)
=4(2-d_K(p))-2(1+d_K(p)).
\]

Every vertex of $K$ belongs to $P$ or $Q$. Summing the two estimates
therefore gives a term $4(2-d_K(w))$ for every $w\in V(K)$, with
$2(1+d_K(p))$ subtracted for each $p\in P\cap V(K)$. The median
$m$ belongs to $P\cap V(K)$ and satisfies $d_K(m)\ge3$, so
\[
\sum_{p\in P\cap V(K)}2(1+d_K(p))
\ge 2(1+d_K(m))\ge8.
\]
It follows that
\begin{equation}\label{eq:component}
\begin{aligned}
&4\sum_{q\in Q\cap V(K)}(2-d_F(q))
 -6\sum_{p\in P\cap V(K)}(d_F(p)-1)\\
&\qquad\le 4\sum_{w\in V(K)}(2-d_K(w))-8=0.
\end{aligned}
\end{equation}
The last equality holds because $K$ is a tree, so
$\sum_{w\in V(K)}d_K(w)=2(|V(K)|-1)$ and hence
$\sum_{w\in V(K)}(2-d_K(w))=2$.

Summing~\eqref{eq:component} over the components of $J$ shows that
the last line of~\eqref{eq:degreebalance} is nonpositive. Thus
\[
12+12bI+3\sigma_W\le\sum_{q\in Q_J}(4s-\sigma(q)),
\]
as required.
\end{proof}

\begin{lemma}\label{lem:branchchoice}
The set $Q_J$ is nonempty. Let $k=\min_{q\in Q_J}\sigma(q)$. Then
\begin{equation}\label{eq:parameters}
0\le k<s\le b-1,\qquad |Q|\ge3,
\end{equation}
and
\begin{equation}\label{eq:excessupper}
4+4bI+\sigma_W\le s|Q|+(s-k)\frac{|Q|}{3}.
\end{equation}
\end{lemma}
\begin{proof}
By~\eqref{eq:sigmanonnegative}, $\sigma(w)\ge0$ for every $w\in W^*$, so we have
$\sigma_W\ge\sum_{q\in Q_J}\sigma(q)$ by~\eqref{eq:sigmatotal}.
Since $I\ge0$, the bound in~\eqref{eq:branchbudget} gives
\[
12+3\sum_{q\in Q_J}\sigma(q)
\le \sum_{q\in Q_J}(4s-\sigma(q)).
\]
Rearranging gives
\[
12\le 4\sum_{q\in Q_J}(s-\sigma(q)).
\]
If $Q_J$ were empty, the right side would be zero. Thus $Q_J$ is
nonempty. Moreover, the sum is positive, so some $q_0\in Q_J$
satisfies $s-\sigma(q_0)>0$. By the definition of $k$, we have
$0\le k\le\sigma(q_0)<s$. Finally, $s\le b-1$ by~\eqref{eq:dbound}.

Since $q_0\in Q_J=Q\cap V(J)$ by~\eqref{eq:branchsets}, the vertex $q_0$ lies on a branch
of some vertex $x\in B$ that is exceptional for $W$. The three
neighbors of $x$ in $W^*$ all belong to $Q$, so $|Q|\ge3$.
This proves~\eqref{eq:parameters}.

By the definition of $k$, every $q\in Q_J$ satisfies
$4s-\sigma(q)\le4s-k$. Since $4s-k>0$ by~\eqref{eq:parameters}
and $|Q_J|\le|Q|$, we obtain from~\eqref{eq:branchbudget}
\[
12+12bI+3\sigma_W
\le \sum_{q\in Q_J}(4s-\sigma(q))
\le (4s-k)|Q_J|
\le (4s-k)|Q|.
\]
Dividing by three and using $4s-k=3s+(s-k)$ gives~\eqref{eq:excessupper}.
\end{proof}

\section{Recoloring two vertices on distinct branches}\label{sec:recolor}

By Lemma~\ref{lem:branchchoice}, $Q_J=Q\cap V(J)$ is nonempty. Since $k$ is
defined as the minimum of $\sigma$ on $Q_J$, choose $h\in Q_J$ with
$\sigma(h)=k$. By~\eqref{eq:branchsets}, $h\in V(J)$, so $h\in S_x$ for some vertex
$x\in B$ exceptional for $W$. Thus $h$ lies on a branch of $x$.
By Lemma~\ref{lem:median}, the median of $x$ belongs to $P$.
Since $h\in Q$, it is distinct from the median. The three branches
meet only at the median, so exactly one of them contains $h$.

Choose one of the two branches of $x$ that do not contain $h$,
and let $q$ be the neighbor of $x$ at the end of this branch.
Since $x$ is exceptional for $W$, $q\in Q$. Also, $q\in S_x\subseteq V(J)$,
so $q\in Q_J$. By the choice of $h$ and the definition of $k$,
\begin{equation}\label{eq:branchvalues}
\sigma(h)=k,\qquad \sigma(q)\ge k.
\end{equation}

The vertices $h$ and $q$ lie on different branches of $x$. Both
belong to $Q$, and the median belongs to $P$. Thus neither vertex
is the median, and Lemma~\ref{lem:median} shows that $x$ is movable
to $W\setminus\{h,q\}$. Put
\begin{equation}\label{eq:H}
H=(W\setminus\{h,q\})\cup\{x\}.
\end{equation}
Then $G[H]$ is a forest. We next obtain a lower bound on the number
of edges from the other color classes to $H$.

By Lemma~\ref{lem:absorb}, we may color $B\setminus\{x\}$ in $b$
forest classes of size $c+1$ and leave $x$ uncolored. Since $q\in Q=W^*\setminus P$ by~\eqref{eq:PQ}, the definition of $P$ implies that
$q$ is movable to some class $Z\in\A\setminus\{W\}$.
By Lemma~\ref{lem:paths},
$Z$ has a path to $C$ avoiding $W$. Adding the arc $W\to Z$ at the
beginning gives a path $P_0$ from $W$ to $C$. Choose $q$ as the
witness for its first arc and fix a witness for every other arc.
Since $W\in\T_U$, the path passes through $U$. Before reaching $U$
it stays in $\T$, since no arc goes from $\T$ to $\OO$. Every class
after $U$ belongs to $\OO$, since it can reach $C$ along the rest
of the path without passing through $U$.

\begin{lemma}\label{lem:omitted}
There is a family $\LL\subseteq\T\setminus\{W\}$ in which every
class $V$ has size at least $|W|$ and every vertex of $V$ has at
least two neighbors in $H$. This family can be chosen to satisfy
\begin{equation}\label{eq:omitted}
|\LL|\ge t-1-k.
\end{equation}
\end{lemma}

\begin{proof}
Call a class in $(\T\setminus\{W\})\cup\{U\}$ a candidate if it
contains at most one neighbor of $h$. Let $m$ be the number of
noncandidates among these $t$ classes.

Since $h\in Q$, we may count one neighbor of $h$ in each of the
$t-1$ classes of $\T\setminus\{W\}$ and two in each of the
$a-t-1$ classes of $\OO$ by~\eqref{eq:classtable}. So far we have
counted $(t-1)+2(a-t-1)$ neighbors, all outside $U$. Each
noncandidate contains at least two neighbors of $h$. If it belongs
to $\T\setminus\{W\}$, only one of these neighbors has been
counted. If it is $U$, none has been counted. Thus we can count
at least one additional neighbor in each of the $m$ noncandidates,
giving
\[
d_{A\setminus W}(h)\ge (t-1)+2(a-t-1)+m.
\]
The second line of~\eqref{eq:sigmadef} now gives
\[
k=\sigma(h)
=d_{A\setminus W}(h)-\bigl((t-1)+2(a-t-1)\bigr)\ge m.
\]
Of the $t$ classes, exactly $m$ are noncandidates, so the remaining
$t-m\ge t-k$ classes are candidates. For any candidate $V$, the
graph $G[V]$ is a forest and $h$ has at most one neighbor in $V$.
Adding $h$ to $V$ therefore creates no cycle, so $h$ is movable to $V$.

Recall that $P_0$ is the directed path from $W$ to $C$ chosen
before the lemma. If any candidates lie on this path, let $Y$ be
the last one encountered from $W$ to $C$, and let $\LL$ consist
of all candidates except $Y$. Otherwise, let $\LL$ consist of all
candidates. Since there are at least $t-k$ candidates and at most
one is omitted, $|\LL|\ge t-k-1$, proving~\eqref{eq:omitted}.
If $U$ is a candidate, then it lies on $P_0$. Every class after
$U$ on $P_0$ belongs to $\OO$, so none is a candidate. Thus $U$
is the last candidate on $P_0$ and is omitted from $\LL$. It follows
that $\LL\subseteq\T\setminus\{W\}$.

Fix $V\in\LL$. Since $V\in\T\setminus\{W\}$ and $C\notin\T$,
the classes $W$, $V$, and $C$ are distinct.

\medskip
\noindent\textbf{Claim 1.}
There is a partition of $V(G)$ into $r$ induced forests in which
the class labeled $W$ has vertex set $H$, the class labeled $V$
contains all vertices of $V$ and one additional vertex, and the
class labeled $C$ has size $c+1$. All other classes have their
original sizes.

\smallskip
\noindent\emph{Proof of Claim 1.}
Suppose first that $V$ lies outside $P_0$. Since $V\in\LL$ and
$\LL$ consists of candidates, $h$ has at most one neighbor in $V$
and is therefore movable to $V$. The switch along $P_0$ leaves
$V$ unchanged and moves $q$ out of $W$ as the witness for the
first arc. Since $q\ne h$, we may perform the switch and move
$h$ from $W$ to $V$ simultaneously. The class labeled $W$ loses
$q$ and $h$, the class labeled $C$ gains one vertex, and the
class labeled $V$ has vertex set $V\cup\{h\}$.

Now suppose that $V$ lies on $P_0$. Since $V\in\LL$, it is a
candidate, so $P_0$ contains at least one candidate. Recall that
$Y$ is the last candidate on $P_0$ from $W$ to $C$ and that we
excluded $Y$ when forming $\LL$. Thus $V$ precedes $Y$, and the
subpaths from $W$ to $V$ and from $Y$ to $C$ have no class in
common. Since $Y$ is a candidate, it contains at most one
neighbor of $h$, so $h$ is movable to $Y$.

We switch along both subpaths and move $h$ from $W$ to $Y$
simultaneously. The switch from $W$ to $V$ uses $q$ as its first
witness, so the class labeled $W$ loses the two distinct vertices
$q$ and $h$, and the class labeled $V$ gains one vertex while
keeping all its original vertices. If $Y\ne C$, the switch from
$Y$ to $C$ moves one vertex out of $Y$ and adds one vertex to $C$.
The incoming vertex $h$ restores the size of $Y$. If $Y=C$,
the second subpath has length zero, and the class labeled $C$
gains $h$ directly.

In both cases, each receiving class gets one vertex that is movable
to its original vertex set. Deleting any outgoing vertex preserves
the resulting forest. Each class other than $W$, $V$, and $C$ either
sends and receives one vertex or is untouched, so its size is unchanged.
Finally, inserting $x$ into $W\setminus\{q,h\}$ gives $H$
by~\eqref{eq:H}, and $G[H]$ is a forest. This proves Claim~1.

The two recolorings in Claim~1 are illustrated in
Figure~\ref{fig:recoloring}.
\begin{figure}[H]
\centering
\begin{tikzpicture}[x=1cm,y=1cm,font=\small,>=stealth,
  class/.style={draw,rounded corners=2pt,minimum width=8mm,
    minimum height=7mm,fill=white},
  move/.style={->,thick},
  change/.style={font=\footnotesize}]
  \node[anchor=west] at (-1.6,2) {(a) $V$ lies outside $P_0$.};
  \node[class] (aW) at (0,0) {$W$};
  \node[class] (aV) at (6,1.1) {$V$};
  \node[class] (aC) at (12,0) {$C$};
  \node (ax) at (-1.4,0) {$x$};
  \draw[move] (ax)--(aW);
  \draw[move] (aW)--node[below=3pt]
    {switch along $P_0$ (first witness $q$)} (aC);
  \draw[move] (aW.north) to[out=65,in=180]
    node[above] {$h$} (aV.west);
  \node[change,below=4pt] at (aW.south) {$-1$};
  \node[change,above=4pt] at (aV.north) {$+1$};
  \node[change,below=4pt] at (aC.south) {$+1$};

  \begin{scope}[yshift=-3.5cm]
    \node[anchor=west] at (-1.6,2) {(b) $V$ lies on $P_0$.};
    \node[class] (bW) at (0,0) {$W$};
    \node[class] (bV) at (4,0) {$V$};
    \node[class] (bY) at (8,0) {$Y$};
    \node[class] (bC) at (12,0) {$C$};
    \node (bx) at (-1.4,0) {$x$};
    \draw[move] (bx)--(bW);
    \draw[move] (bW)--node[above=3pt] {first witness $q$}
      node[below=3pt] {switch} (bV);
    \draw[->,dashed] (bV)--node[below=3pt] {no switch} (bY);
    \draw[move] (bY)--node[below=3pt] {switch} (bC);
    \draw[move] (bW.north) .. controls (1.2,1.65) and (6.8,1.65) ..
      node[above] {$h$} (bY.north);
    \node[change,below=18pt] at (bW.south) {$-1$};
    \node[change,below=18pt] at (bV.south) {$+1$};
    \node[change,below=18pt] at (bY.south) {$0$};
    \node[change,below=18pt] at (bC.south) {$+1$};
  \end{scope}
\end{tikzpicture}
\caption{The two recolorings in Claim~1 of Lemma~\ref{lem:omitted}.
Boxes denote class labels. Horizontal arrows represent subpaths
of $P_0$, with intermediate classes suppressed. Solid horizontal
arrows indicate switches, and the dashed subpath is not switched.
The other arrows show the moves of $h$ and $x$.
The numbers give the changes in class sizes after $x$ is inserted
into $W\setminus\{h,q\}$, producing $H$.
Panel (b) shows $Y\ne C$. If $Y=C$, the last subpath has length
zero and $C$ receives $h$ directly.}
\label{fig:recoloring}
\end{figure}
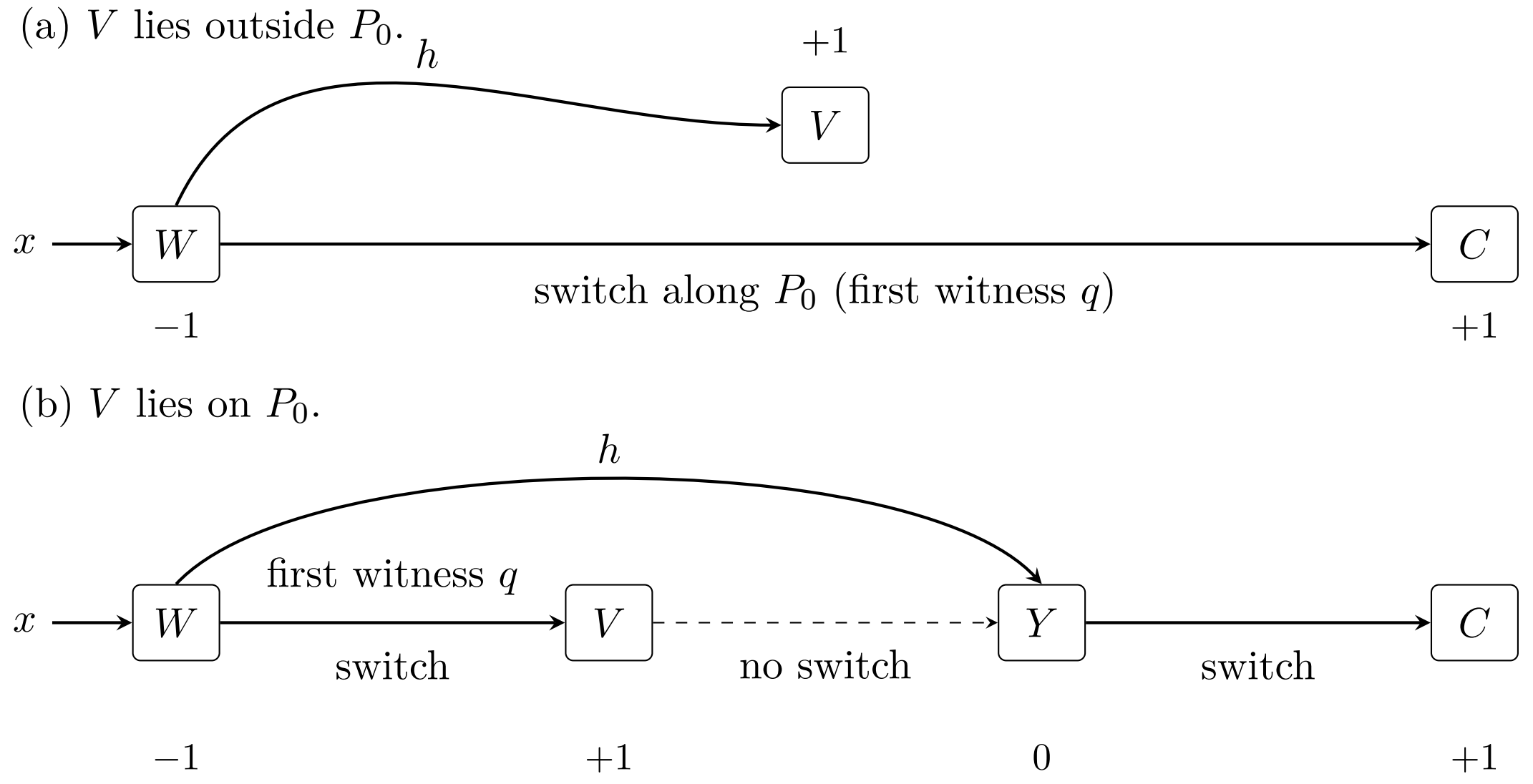

\medskip
\noindent\textbf{Claim 2.}
We have $|V|\ge|W|$ and $d_H(v)\ge2$ for every $v\in V$.

\smallskip
\noindent\emph{Proof of Claim 2.}
In the partition given by Claim~1, the classes labeled $W$, $V$,
and $C$ have sizes $|W|-1$, $|V|+1$, and $c+1$, respectively.
If $|V|<|W|$, then $|V|=c$ and $|W|=c+1$ by~\eqref{eq:basic}.
The three new sizes are then $c$, $c+1$, and $c+1$, so the
partition is an equitable tree-$r$-coloring of $G$. This
contradiction gives $|V|\ge|W|$.

Suppose that some $v\in V$ has at most one neighbor in $H$.
By Claim~1, the class labeled $V$ still contains $v$. Moving $v$
from that class to $H$ preserves both forests and restores the
original sizes of the classes labeled $V$ and $W$. The class
labeled $C$ has size $c+1$, and every other class now has its
original size. This again gives an equitable tree-$r$-coloring
of $G$. Therefore $d_H(v)\ge2$ for every $v\in V$. This proves
Claim~2.

Since $V\in\LL$ was arbitrary, Claim~2 shows that every class
$V\in\LL$ has size at least $|W|$ and every vertex of $V$ has at
least two neighbors in $H$.
\end{proof}

\section{The final degree count}\label{sec:final}

Let $\LL$ be a family of classes given by Lemma~\ref{lem:omitted},
and put
\begin{equation}\label{eq:Ldef}
L=\bigcup_{V\in\LL}V,\qquad \ell=|\LL|.
\end{equation}
The classes in
$\LL$ are pairwise disjoint and each has size at least $|W|$.
Hence
\[
|L|=\sum_{V\in\LL}|V|\ge\ell|W|.
\]

Since $\LL\subseteq\T\setminus\{W\}$, we have $L\subseteq A\setminus W$.
By~\eqref{eq:H}, $H\subseteq W\cup\{x\}$ with $x\in B$, so $H$
is disjoint from $L$. Every vertex of $L$ belongs to a class in
$\LL$ and therefore has at least two neighbors in $H$ by
Lemma~\ref{lem:omitted}. Summing over $L$ gives
\[
e(L,H)=\sum_{v\in L}d_H(v)\ge2|L|.
\]

Put
\begin{equation}\label{eq:Wparts}
W'=W^*\setminus\{h,q\},\qquad Z=W\setminus W^*.
\end{equation}
Since
$h,q\in W^*$, the set $H$ is the disjoint union of $W'$, $Z$,
and $\{x\}$. Counting the edges from $L$ to these three sets gives
\begin{equation}\label{eq:demand}
2\ell|W|\le 2|L|\le e(L,H)=e(L,W')+e(L,Z)+d_L(x).
\end{equation}

To estimate $e(L,W')$, fix $w\in W'$. Since
$\LL\subseteq\T\setminus\{W\}$, we have $L\subseteq A\setminus W$.
We bound $d_L(w)$ by counting the neighbors of $w$ in $A\setminus W$
and subtracting a lower bound on how many of them lie outside $L$.

If $w\in P$, then $d_{A\setminus W}(w)=2(a-1)+\sigma(w)$
by~\eqref{eq:sigmadef}. There are $a$ accessible classes.
Excluding $W$ and the $\ell$ classes in $\LL$ leaves
$a-1-\ell$ classes, all disjoint from $L$. By the first row
of~\eqref{eq:classtable}, $w$ has at least two neighbors in each
of these classes. Thus at least $2(a-1-\ell)$ of its neighbors
in $A\setminus W$ lie outside $L$.

If $w\in Q$, then~\eqref{eq:sigmadef} gives
$d_{A\setminus W}(w)=(t-1)+2(a-t-1)+\sigma(w)$.
The family $\T$ has $t$ classes. Excluding $W$ and the $\ell$
classes in $\LL$ leaves $t-1-\ell$ classes in $\T$, all disjoint
from $L$. By the second row of~\eqref{eq:classtable}, $w$ has at
least one neighbor in each of these classes.
By~\eqref{eq:families}, the family $\OO$ consists of the $a-t-1$
accessible classes outside $\T\cup\{U\}$. These classes are also
disjoint from $L$, and the same row gives at least two neighbors
of $w$ in each of them. Thus at least $(t-1-\ell)+2(a-t-1)$ of
the neighbors of $w$ in $A\setminus W$ lie outside $L$.

Subtracting the lower bound for neighbors outside $L$ from the
total in $A\setminus W$ gives, for $w\in P$,
\[
\begin{aligned}
d_L(w)
&\le 2(a-1)+\sigma(w)-2(a-1-\ell)\\
&=2\ell+\sigma(w).
\end{aligned}
\]
For $w\in Q\setminus\{h,q\}$, the subtraction gives
\[
\begin{aligned}
d_L(w)
&\le (t-1)+2(a-t-1)+\sigma(w)
     -\bigl((t-1-\ell)+2(a-t-1)\bigr)\\
&=\ell+\sigma(w).
\end{aligned}
\]

By~\eqref{eq:PQ}, $W^*$ is the disjoint union of $P$ and $Q$.
Since $h,q$ are distinct vertices of $Q$, \eqref{eq:Wparts} gives
$W'=P\cup(Q\setminus\{h,q\})$. Thus $W'$ contains all of $P$
and exactly $|Q|-2$ vertices of $Q$.
By~\eqref{eq:sigmatotal}, the sum of $\sigma(w)$ over $W'$ is
$\sigma_W-\sigma(h)-\sigma(q)$. By~\eqref{eq:branchvalues}, $\sigma(h)=k$ and $\sigma(q)\ge k$.
Since each edge between $L$ and $W'$ has exactly one endpoint in
$W'$, we have $e(L,W')=\sum_{w\in W'}d_L(w)$. Splitting this sum
over $P$ and $Q\setminus\{h,q\}$ and applying the bound on $d_L(w)$
for each set gives
\begin{equation}\label{eq:Lactive}
\begin{aligned}
e(L,W')
&= \sum_{w\in P}d_L(w)
   +\sum_{w\in Q\setminus\{h,q\}}d_L(w)\\
&\le \sum_{w\in P}\bigl(2\ell+\sigma(w)\bigr)
   +\sum_{w\in Q\setminus\{h,q\}}\bigl(\ell+\sigma(w)\bigr)\\
&= 2\ell|P|+\ell(|Q|-2)+\sum_{w\in W'}\sigma(w)\\
&= \ell(2|P|+|Q|-2)+\sigma_W-\sigma(h)-\sigma(q)\\
&\le \ell(2|P|+|Q|-2)+\sigma_W-2k.
\end{aligned}
\end{equation}

Substituting~\eqref{eq:Lactive} into~\eqref{eq:demand} gives
\[
2\ell|W|
\le \ell(2|P|+|Q|-2)+\sigma_W-2k+e(L,Z)+d_L(x).
\]
By~\eqref{eq:PQ}, $W^*$ is the disjoint union of $P$ and $Q$.
By~\eqref{eq:Wparts} and~\eqref{eq:isolated}, the remaining
vertices of $W$ form the set $Z=W\setminus W^*$ of size $I$.
Hence $|W|=|P|+|Q|+I$.

Subtracting $\ell(2|P|+|Q|-2)$ from both sides of the preceding
inequality, we obtain on the left
\[
\begin{aligned}
2\ell|W|-\ell(2|P|+|Q|-2)
&=2\ell(|P|+|Q|+I)-\ell(2|P|+|Q|-2)\\
&=\ell|Q|+2\ell I+2\ell.
\end{aligned}
\]
Thus
\[
\ell|Q|+2\ell I+2\ell
\le \sigma_W-2k+e(L,Z)+d_L(x).
\]
Adding $2k$ to both sides and subtracting $2\ell I+2\ell$ gives
\begin{equation}\label{eq:finalbalance}
\ell|Q|+2k
\le \sigma_W+\bigl(e(L,Z)-2\ell I\bigr)
             +\bigl(d_L(x)-2\ell\bigr).
\end{equation}

Recall from~\eqref{eq:Wparts} that $Z=W\setminus W^*$ is the set
of isolated vertices of $G[W]$. Since $W\in\T$,
Lemma~\ref{lem:paths} excludes an arc from $W$ to any class in
$\OO$. Fix $z\in Z$ and let $V\in\OO$. If $z$ had at most one
neighbor in $V$, adding $z$ to $V$ would preserve the forest.
Since $z\in W$, it would then witness an arc from $W$ to $V$.
Thus $z$ has at least two neighbors in every class of $\OO$.

By~\eqref{eq:families}, $\OO$ has $a-t-1$ classes, all outside
$\T$. Since $\LL\subseteq\T$, these classes are disjoint from $L$.
Hence $z$ has at least $2(a-t-1)$ neighbors outside $L$, giving
\[
d_L(z)\le d_G(z)-2(a-t-1)\le 2r-1-2(a-t-1).
\]
By~\eqref{eq:isolated} and~\eqref{eq:Wparts}, $|Z|=I$. Summing
over $Z$ therefore gives
\[
e(L,Z)=\sum_{z\in Z}d_L(z)\le \bigl(2r-1-2(a-t-1)\bigr)I.
\]
Subtracting $2\ell I$ and using $r=a+b$ by~\eqref{eq:partition},
we obtain
\begin{equation}\label{eq:Zcount}
\begin{aligned}
e(L,Z)-2\ell I
&\le \bigl(2r-1-2(a-t-1)-2\ell\bigr)I\\
&=\bigl(2a+2b-1-2a+2t+2-2\ell\bigr)I\\
&=\bigl(2b+2(t-\ell)+1\bigr)I.
\end{aligned}
\end{equation}

Since $\ell=|\LL|$, inequality~\eqref{eq:omitted} gives
$\ell\ge t-1-k$, or equivalently $t-\ell\le1+k$.
By~\eqref{eq:parameters}, $k<s\le b-1$. Since $k$ and $b$ are
integers, $k\le b-2$. Thus the coefficient of $I$
in~\eqref{eq:Zcount} satisfies
\[
2b+2(t-\ell)+1
\le 2b+2(1+k)+1
\le 2b+2(b-1)+1
=4b-1.
\]
Since $I\ge0$, it follows that
\begin{equation}\label{eq:Zupper}
e(L,Z)-2\ell I\le 4bI.
\end{equation}

The vertex $x$ has exactly three neighbors in $W^*$ because it
is exceptional for $W$. By~\eqref{eq:basic}, it also has at least
two neighbors in each of the $a-1-\ell$ accessible classes outside
$\LL\cup\{W\}$. All these neighbors lie outside $L$. Subtracting
their number from the degree bound $2r-1$ gives
\begin{equation}\label{eq:xcount}
d_L(x)\le2r-1-3-2(a-1-\ell)=2\ell+2b-2.
\end{equation}

Using~\eqref{eq:Zupper} and~\eqref{eq:xcount}
in~\eqref{eq:finalbalance}, we get
\[
\ell|Q|+2k\le\sigma_W+4bI+2b-2.
\]
By~\eqref{eq:omitted}, $\ell\ge t-1-k$. Since $s=t-b$
by~\eqref{eq:sigmatotal}, this is $\ell\ge s+(b-1-k)$.
Using this lower bound on $\ell$ and writing
$2k-2b+2=-2(b-1-k)$ gives
\begin{equation}\label{eq:excesslower}
\begin{aligned}
4bI+\sigma_W
&\ge \ell|Q|-2(b-1-k)\\
&\ge s|Q|+(b-1-k)(|Q|-2).
\end{aligned}
\end{equation}

By~\eqref{eq:parameters}, we have $b-1-k\ge s-k>0$.
The same estimate gives $|Q|\ge3$, so
$|Q|-2\ge|Q|/3$. Hence~\eqref{eq:excesslower} yields
$4bI+\sigma_W\ge s|Q|+(s-k)|Q|/3$. Combining this lower
bound with the upper bound in~\eqref{eq:excessupper} gives
\[
4+4bI+\sigma_W
\le s|Q|+(s-k)\frac{|Q|}{3}
\le 4bI+\sigma_W.
\]
This contradiction completes the induction on $r$ and proves
Theorem~\ref{thm:main}.\hfill$\square$

\end{document}